\documentclass[11pt,a4paper]{article}
 \usepackage{euscript}

\usepackage{amsmath}
\usepackage{graphicx}
\usepackage{amsthm}
\usepackage{amssymb}
\usepackage{epsfig}

\usepackage{url}

\theoremstyle{theorem}
\newtheorem{theorem}{Theorem}[section]
\newtheorem{corollary}[theorem]{Corollary}

\newtheorem{lemma}[theorem]{Lemma}
\newtheorem{proposition}[theorem]{Proposition}

\newtheorem{definition}{Definition}[section]

\newtheorem{example}{Example}[section]

\numberwithin{equation}{section}

\begin{document}
 \title{On the Field of Rational Pseudo-differential Operators}
 \author{Masood Aryapoor}
  \maketitle
   \begin{abstract}
\noindent
In this paper we study  some properties of the field of rational pseudo-differential operators on a field and some other related rings. As an application
we reconstruct the Kac co-cycle on the Lie algebra of differential operators on a circle. 
\end{abstract}

                                                              \begin{section}{Introduction}
  
 Suppose that $F$ is a differential field of characteristic zero. One can associate three noncommutative rings to $F$. They are\\ 
 (1) the ring of differential operators, denoted by $F\langle \partial \rangle$,\\
 (2) the field of rational pseudo-differential operators, denoted by  $F(\partial)$,\\
 (3)  the field of formal pseudo-differential operators, denoted by $F((\partial^{-1}))$.\\
 The goal of this paper is to study these objects as noncommutative counterparts of \\
 (1) $F[x]$, the ring of polynomials over $F$,\\
 (2) $F(x)$, the field of rational functions over $F$,\\
 (3) $F((x^{-1}))$, the field of Laurent series over $F$.\\
 In particular we define the analogue of zeros and poles for rational pseudo-differential operators( section 2.3.).
 
 In the case where $F$ is linearly differentially closed we obtain a decomposition 
 $$F(\partial)=F\langle \partial \rangle +F\langle \partial^{-1} \rangle $$ 
 where $F\langle \partial^{-1} \rangle$ is the ring generated by $F$ and $\partial^{-1}$ (Theorem 3.8.). 
 
 In section 5, we show that how one can embed the field of rational pseudo-differential operators into the Calkin algebra. This, in particular, gives us a Lie algebra 2 co-cycle on 
 $F(\partial)$. Furthermore we show that this 2 co-cycle gives us a  2 co-cycle on the Lie algebra of vector fields on the circle which coincides with the Kac co-cycle.
 
 \bf{Acknowledgement}\\
 I would like to thank my advisor, professor Kapranov, for his support and useful discussions. Also this work was partially done during my stay at the
  Max-Planck-Insitut f\"{u}r Mathematik and this reasearch was partially  supported by the Max-Planck-Insitut f\"{u}r Mathematik.   
  
 \end{section}

              \begin{section} {Differential operators}
   
                                                    \begin{subsection}{Ring of differential operators}

Suppose that $F$ is a differential field of characteristic zero. By a differential field we mean a field with a derivation. We denote the derivative of $a\in F$ by $a'$.
Suppose that $C$ is the set of elements of $F$ whose derivatives are zero.
It is easy to see that $C$ is a filed, called the field of constants. We assume that $C$ is algebraically closed and $C\neq F$
(which implies that $dim_C(F)=\infty$). For more on differential algebras see  [Kaplansky].

The ring of differential operators over $F$ is defined to be the following ring,\\
$$ F\langle \partial \rangle = \langle F, \partial ;\partial a-a \partial =a' \quad \text{for every} \quad a\in F \rangle.$$
Elements of $F\langle \partial \rangle$ are called differential operators.\\
Every differential operator $P\in F\langle \partial \rangle  $ has a unique presentation as $P=a_0+a_1\partial +
\cdots+a_n \partial^n$ where $a_0,...,a_n\in F\langle \partial \rangle$ and $a_n\neq 0$. We call $n$ the degree of $P$ and denote it by $deg(P)$. It is straightforward
 to check that\\
 1)For every $P,Q\in F\langle \partial \rangle $, $deg(PQ)=deg(P)+deg(Q)$.\\
 2)For every  $P,Q\in F\langle \partial \rangle $, $deg(P+Q) \leq \max (deg(P),deg(Q))$.\\
Using the degree map, one can see that $F\langle \partial \rangle $ is a filtered ring. The filtration is given by
$F\langle \partial \rangle_{\leq n}=\{P\in F\langle \partial \rangle| deg(P)\leq n\}$. It is easy to see that $gr(F\langle \partial \rangle)=F[x]$.\\

We summarize some of the well-known properties  of  $F\langle \partial \rangle$ in the following proposition,
\begin{proposition}
For and differential field $F$ we have,\\
(a) $F\langle \partial \rangle $ is a domain, i.e. it does not have any zero divisors.\\
(b) Every left(right) ideal of $F\langle \partial \rangle$ is principal.\\
(c) $F\langle \partial \rangle$  is a simple $C$-algebra.\\
(e) For any $a\in F\backslash C$,  the centralizer of $a$ in $F\langle \partial \rangle$ is $F$. 
\end{proposition}
We can consider a differential operator $P$ as a function $P:F\to F$. More precisely,  if $P=a_0+a_1\partial +
\cdots+a_n \partial^n$ and $a\in F$ we define, $$P(a)=a_0a+a_1a'+\cdots+a_na^{(n)}.$$ 
An element $a\in F$ is called a zero of $P$ if $P(a)=0$. The set of zeros of $P$ is denoted by $Zer(P)$.
It is well-known that,
\begin{proposition}
For any differential operator $P$, $Zer(P)$  is a vector space over $C$ of dimension $\leq deg(P)$.
\end{proposition} 
See [Kaplansky] for a proof.

It is also known that,
\begin{proposition}
Given any  $C$-vector subspace $V\subset F$ of dimension $n$, there is a unique differential operator $P_V=a_0+a_1\partial +
\cdots+a_{n-1}\partial^{n-1}+\partial^{n}$ such that $Zer(P_V)=V$. 
\end{proposition}
 
\end{subsection}

                                                                            \begin{subsection}{Wronskian}

Given $x_1,\dots,x_n \in F$, the wronskian of them, denoted by $w(x_1,...,x_n)$, is defined to be the determinant of the following matrix,\\

\(
\begin{bmatrix}
 x_1& \cdots & x_n \\
 x_1'&\cdots & x_n'\\
 \vdots&\vdots& \vdots \\
 x_1^{(n-1)}&\cdots & x_n^{(n-1)}
 \end{bmatrix}
 \) 
\newline 
where $x^{(m)}=(x^{(m-1)})'$ and $x^{(0)}=x$.\\
Wronskian is a useful tool in studying differential rings( See [Van der Put] for more on wronskian). As an example we recall the following,
\begin{proposition}
Suppose that $x_1,\dots,x_n \in F$. Then $x_1,\dots,x_n \in F$ are linearly dependent over $C$ if and only if $w(x_1,...,x_n)=0$.
\end{proposition}
This enables us to find $P_V$ for a finite dimensional vector space $V\subset F$ over $C$(proposition 2.3.).
\begin{proposition}
Suppose that   $x_1,...,x_n\in F$ is a basis for the finite dimensional $C$-vector space $V\subset F$. Suppose that \\
\(
\begin{bmatrix}
y& x_1& \cdots & x_n \\
y'& x_1'&\cdots & x_n'\\
\vdots & \vdots&\vdots& \vdots \\
y^{(n)} & x_1^{(n)}&\cdots & x_n^{(n)}
 \end{bmatrix}
 \)= $a_{n}y^{(n)}+\cdots+a_1y'+a_0y.$\\
 Then $Z(a_{n}\partial^{n}+\cdots+a_1\partial+a_0)=V$.
\end{proposition}

                                                     \begin{subsection}{Linearly Differentially closed fields}
						     
We start with  the following definition,
\begin{definition}
A differential field $F$ is called linearly differentially closed if any differential operator $P\in F\langle \partial \rangle $ has a nonzero solution in $F$.
\end{definition}
>From Galois theory of differential equations, one can see that every differential field can be embedded into a linearly differentially closed field(See [Van
der Put] or [Magid]).\\
In this section we suppose that $F$ is linearly differentially closed. First of all, it is easy to see that\\ 
(1) every differential operator can be written as the product of degree one differential operators,\\
(2) every nonzero differential operator $P:F\to F$ is onto,\\
(3) for any nonzero differential operator $P$, we have $dim_C(Z(P))=deg(P)$.\\
\begin{definition}
Suppose that $P=a_0+a_1\partial +
\cdots+a_{n-1}\partial^{n-1}+a_n\partial^{n}$ is a differential operator. A decomposition of $P$ into irreducible factors is a decompostion of $P$ of the form
$P=a_n(\partial+x_1)\cdots(\partial+x_n)$.
\end{definition}
\begin{proposition}
For any differential operator $P$ of degree $n$, there is a one-to-one correspondence between decompositions of $P$ into irreducible factors 
and the points of $Fl(1,2,\cdots,n;Z(P))$, the full flag variety of $Z(P)$ over $C$.
\end{proposition}  
\begin{proof}
The correspondence is given by sending a decomposition $$P=a_n(\partial+x_1)\cdots(\partial+x_n)$$ into $$Z(\partial+x_n)\subset Z((\partial+x_{n-1})(\partial+x_n))\subset 
\cdots \subset Z((\partial+x_2)\dots(\partial+x_n))\subset Z(P).$$
Proposition 2.3. shows that this is a one-to-one correspondence.
\end{proof}

\end{subsection}

\end{subsection}

\end{section}

                              \begin{section}{Field of rational pseudo-differential operators}

                                                          \begin{subsection}{Ore property}
We recall the definition of Ore domains,
\begin{definition}
A domain $R$ is called  an Ore domain  if for any $r,s\in R\backslash \{0\}$, $Rr\bigcap Rs \neq \{ 0 \}$ and $rR\bigcap sR\neq \{ 0\}$.
\end{definition}
It is a well-known fact that every Ore domain $R$  has a filed of quotients $Q(R)$. Elements of $Q(R)$ can be considered as $r^{-1}s$ where $r,s\in R$. For more on Ore property
see [Lam].
\end{subsection}

                                                         \begin{subsection}{Rational pseudo-differential operators}
 It is well-known that,
 \begin{proposition}
 The ring of differential operators over $F$ is an Ore domain.
 \end{proposition}
 See [Bjork].\\
 The quotient field of $F\langle \partial \rangle$ is called the field of rational pseudo-differential operators and denoted by $F(\partial)$. Its elements are called rational
 pseudo-differential operators. Every  rational pseudo-differential operator can be presented as $P^{-1}Q$ where $P,Q$ are differential operators. 
 \begin{definition}
 A minimal presentation of $f\in F(\partial)$ is a presentation $f=P^{-1}Q$ where $P,Q$ are differential operators and $P$ has the smallest possible degree.
 \end{definition}  
It is easy to see that,
\begin{lemma}
 If $P^{-1}Q$ and $P_1^{-1}Q_1^{-1}$ are two minimal presentation of $f\in F(\partial)$, then there is $a\in F$ such that  $P_1=aP$ and $Q_1=aQ$.
\end{lemma}
We call the degree of $P$ the length of $f$  and denote it by $L(f)$ if $P^{-1}Q$ is a minimal presentation of $f$.

 \end{subsection}

	                                                \begin{subsection}{Zeros and Poles}
							   
A rational pseudo-differential operator cannot be considered as a function on $F$. Nevertheless we
can define poles and zeros of a rational pseudo-differential operators as follows.
 \begin{definition}
 For $f\in F(\partial)$, we define $Zer(f)=Zer(Q)$ and $Pol(f)=Zer(P)$ where $f=P^{-1}Q$ is a minimal presentation of $f$.
 \end{definition}
 Note that this definition makes sense by lemma 3.2.
 Next we show that $Zer(f)$ and $Pol(f)$ can be completely arbitrary.
 We begin with a lemma,
\begin{lemma}

 Let $V$ and $W$ be two finite dimensional $C$-vector subspaces of $F$. Then there is some $a\in F\backslash\{0\}$ such that $aV\bigcap W=\{0\}$.

\end{lemma}

\begin{proof} 

Let $Z$ be a $C$-vector subspace of $F$ such that $W\bigoplus Z=F$.  Let $a_1,...,a_n$ be a basis for $V$. The vector spaces $a_1^{-1}Z,...,a_n^{-1}Z$ have 
finite co-dimensions. This  implies that  $a_1^{-1}Z\bigcap a_2^{-1}Z\bigcap...\bigcap a_n^{-1}Z$ is a nonzero vector subspace. Therefore for 
$0\neq a\in a_1^{-1}Z\bigcap a_2^{-1}Z\bigcap...\bigcap a_n^{-1}Z$, we have $aV\subset Z$, which implies that $aV\bigcap W=\{0\}$.

\end{proof}
Using this lemma we have,
\begin{proposition}

 Given two finite dimensional $C$-vector subspaces $V$ and $W$ of $F$, there is $f\in F(\partial)$ with $Pol(f)=V$ and $Zer(f)=W$.\\

\end{proposition}

\begin{proof}
By the above lemma, there is $a\in F\backslash \{0\}$ such that $f=P_{V}^{-1}(aP_W)$ is a minimal presentation. So
$Pol(f)=V$ and $Zer(f)=W$.
\end{proof}

\begin{example}
There is a rational pseudo-differential operator $f$ such that $Zer(f)=Pol(f)$. It has the form $P^{-1}aP$ where $P$ is an appropriate differential operator and $a\in F$.
\end{example}

 \end{subsection}
 
                                                \begin{subsection}{The ring of integration operators}
						
 One can think of $\partial^{-1}$ as an integration. This leads to the following definition,
 \begin{definition}
 The ring of integration operators over $F$ is defined to be the subring of $F(\partial)$ generated by $F$ and $\partial^{-1}$ and denoted by $F\langle \partial^{-1} \rangle$.
 \end{definition}
 Elements of $F\langle \partial^{-1} \rangle$ are called integration operators. We want to explain the structure of $F\langle \partial^{-1} \rangle$.
 \begin{proposition}
Assume that the derviation of $F$ is onto. Then the map $\phi: F\bigoplus (F\otimes_C F)\to F\langle \partial^{-1} \rangle$, 
defined by $$\phi(a+\sum_{i=1}^{n}a_i\otimes b_i )=  a+\sum_{i=1}^{n}{a_i\partial^{-1} b_i}$$ is a bijection.
\end{proposition}
\begin{proof}
Thanks to the identity $\partial^{-1}a'\partial^{-1}=a\partial^{-1}-\partial^{-1}a$, $\phi$ is onto. Suppose that $\phi(a+\sum_{i=1}^{n}{a_i\otimes b_i})=0$. We can assume that
$b_1,...,b_n$ are $C$-linearly independent. It is easy to see that $a=0$ and $\sum_{i=1}^{n}{a_ib^{(k)}_i}=0$ for any $i=0,1,...$. But this implies that $a_1=\dots=a_n=0$,
since $b_1,...,b_n$ are $C$-linearly independent( proposition 2.4.).
\end{proof}
This proposition implies that
\begin{corollary}
Assume that the derviation of $F$ is onto. Then $F\langle \partial^{-1} \rangle$ has the following presentation over $F$,
$$ F\langle \partial^{-1} \rangle= \langle F, \partial^{-1} ;a\partial^{-1} - \partial^{-1}a =\partial^{-1}a'\partial^{-1} \quad \text{for every} \quad a\in F\rangle.$$
\end{corollary}

\begin{proposition}
(a) If $I=a+\sum_{i=1}^{n}{a_i\partial^{-1} b_i}\in F\langle \partial^{-1} \rangle$ then $L(I)\leq n$. Moreover there is a presentation of $I$ with  $L(I)$ tensors.\\
(b) Let $I=a+\sum_{i=1}^{n}{a_i\partial^{-1} b_i}\in F\langle \partial^{-1} \rangle$ with $n=L(I)$. Then $Pol(I)=Ca_1+\cdots + Ca_n$.
\end{proposition}

\begin{proof}
 (a) It is easy to see.\\
 (b) If $I=P^{-1}Q$, then $PI=Q$. This implies that $\sum_{i=1}^{n}P(a_i)\partial^{-1}b_i=0$, because
$$PI=Pa+\sum_{i=1}^{n}(Pa_i-P(a_i))\partial^{-1}b_i+\sum_{i=1}^{n}P(a_i)\partial^{-1}b_i$$
and $(Pa_i-P(a_i))\partial^{-1}b_i\in F\langle \partial \rangle $ for each $i$.\\
 Since $b_1,...,b_n$ are $C$-linearly independent, we have $P(a_1)= \cdots=P(a_n)=0$. Conversely, if   $P(a_1)= \cdots=P(a_n)=0$ then $PI\in F\langle \partial \rangle $.
This finishes the proof.
 \end{proof}

\end{subsection}

                                                        \begin{subsection}{Decomposition Theorem}
							
Suppse that $F$ is linearly differentially closed. Then it turns out that $F(\partial)$ has a simple description,
\begin{theorem}
Suppose that $F$ is linearly differentially closed. Then $$F(\partial)=F\langle \partial \rangle +F\langle \partial^{-1} \rangle $$
\end{theorem}
\begin{proof}
Clearly $F(\partial)=F\langle \partial \rangle+F(\partial)_{\leq 0}$. So it is enough to show that $F(\partial)_{\leq 0}=F\langle \partial^{-1} \rangle $.
This follows from the fact that every differential operator can be written as a product of degree one differential operators and the identity 
$(\partial^{-1}a'-a)^{-1}=\partial^{-1}(a^{-1})'-a^{-1}$, $a\in F$.
\end{proof}

\end{subsection}

\end{section}

         \begin{section}{Field of formal pseudo-differential operators}

                                                            \begin{subsection}{Valuation subrings of $F(\partial)$}

 We recall the definition of valuation subrings,
 \begin{definition}
 Suppose that $E$ is  a field( not necessarily commutative). A valuation subring of $E$ is a subring $R$ of $E$ which is invariant under inner automorphisms and for any 
 $0\neq x\in E$, we have $x\in R$ or $x^{-1}\in R$. 
 \end{definition}
 It is easy to see that the set of rational pseudo-differential operators of nonpositive degree, denoted by $F(\partial)_{\leq 0}$,
  is a valuation subring of $F(\partial)$ containing $F$. In fact we have,
\begin{proposition}
The only nontrivial valuation subring of $F(\partial)$ containing $F$ is $F(\partial)_{\leq 0}$.
\end{proposition}
\begin{proof}
Because of the identity $(\partial^{-1}a'-a)^{-1}=\partial^{-1}(a^{-1})'-a^{-1}$, for any valuation subring $F\subset R$ of $F(\partial)$ we have $\partial^{-1} \in R$. 
It is now easy to see that $F(\partial)_{\leq 0}\subset R$. 
\end{proof}
 
\end{subsection}
                        \begin{subsection}{Completion of the field of rational pseudo-differential operators}

It is well-known that $Ord:F(\partial)\to \mathbb{Z}\bigcup \{ \infty \}$, defined by $Ord(P^{-1}Q)=deg(Q)-deg(P)$($P,Q\in F\langle \partial  \rangle$),
is a valuation on $F(\partial)$. The completion of $F(\partial)$ using this valuation is called the field of formal pseudo-differential operators and denoted by
$F((\partial^{-1}))$. Its elements are called formal pseudo-differential operators.

One can construct $F((\partial^{-1}))$ more concretely as follows.
As a set $F((\partial^{-1}))$ is the set of formal sums $\sum_{i=n}^{\infty}{a_i\partial^{-i}}$ where $n\in \mathbb{Z}$ and $a_i\in F$. The multiplication is given by,
$$ \partial^{n}a=\sum_{i=0}^{\infty}{n\choose_{i}}  a^{(i)} \partial^{n-i}$$
where $n\in \mathbb{Z}$ and $a\in F$.

For $f=\sum_{i=n}^{\infty}{a_i\partial^{-i}}\in F((\partial^{-1}))$,  we define $D(f)=-\sum_{i=n}^{\infty}{ia_i\partial^{-i-1}}$. It is easy to see that $D$ is a derivation of
$F((\partial^{-1}))$. One can check that,
\begin{proposition}
For any $f\in F((\partial^{-1}))$ and $a\in F$ we have,
$$fa= \sum_{i=0}^{\infty}{\frac{1}{i!}a^{(i)}D^i(f)}$$
\end{proposition} 
 
\end{subsection}

\end{section}

                                            \begin{section}{Embedding of the field of rational pseudo-differential  operators into the Calkin algebra}

                                                           \begin{subsection}{Calkin Algebra}
							   
As we mentioned a rational pseudo-differential operator cannot be considered as a function on $F$. However they can be considered as multi-valued functions. 
In this section we use this idea( see [Wells] for psuedo-differential operators in the analytic context).\\							   
In this section, we assume that $C$ is an arbitrary field and $V$ is a vector space
of infinite dimension over  $C$. Set $End_C(V)$ to be the set of  $C$-linear maps form $V$ to
itself. Suppose that $I(V)$ is the set of $C$-linear maps $L\in End_C(V)$  having finite dimensional images. It is easy to see that,
\begin{proposition}
The ideal $I(V)$ is an (two-sided)ideal of $End_C(V)$.
\end{proposition} 
  \begin{definition}
  The Calkin algebra of  $V$ over $C$ is defined to be $M_C(V)=End_C(V)/I(V)$.
  \end{definition}
  It is easy to see that,
 \begin{proposition}
 The element $L+I(V)\in M_C(V)$($L\in End_C(V)$) is invertible if and only if the kernel and co-kernel of $L$ are finite dimensional.
 \end{proposition}   
  
 \end{subsection}
 
                     \begin{subsection}{Embedding of $F(\partial)$ into the Calkin algebra}
		     
In this section we assume that $F$ is linearly differentially closed.  Clearly we have a homomorphism $\theta: F\langle \partial \rangle \to M_C(F)$. 
\begin{proposition}
The homomorphism $\theta$ can be extended to obtain an embedding $F(\partial)\to M_C(F)$, also denoted by $\theta$.
\end{proposition} 
\begin{proof}
This easily follows from proposition 5.2.
\end{proof}
 One can ask whether it is possible to extend this embedding to an embedding of $F((\partial^{-1}))$ into the Calkin algebra of $F$. 
In general, since an element of $F((\partial^{-1}))$ has infinitely many terms,
it seems impossible to do that. However if $F$ is occupied with some topology it might be possible. We explain one situation in with it is possible to do so.

Consider $\mathbb{C}((x))$ as 
a differential field where the derivative is $\frac{\partial}{\partial x}$. Then
\begin{proposition}
There is an embedding 
$$\theta: \mathbb{C}((x))((\partial^{-1}))\to M_{\mathbb{C}}(\mathbb{C}[[x]])$$
extending the natural embedding
$$\theta: \mathbb{C}[[x]]\langle \partial,\partial^{-1}\rangle\to M_{\mathbb{C}}(\mathbb{C}[[x]]).$$
\end{proposition}
\begin{proof}
Roughly speaking, the reason is that the integration increases the degree of elements in $\mathbb{C}[[x]]$ where $deg(\sum_{i=n}^{\infty}a_nx^{n})=n$($a_n\neq 0$). 
More precisely one can see that 
$$\theta(f)(A)=\sum_{k=0}^{n}\frac{\partial^{k}A}{\partial x^{k}}+\sum_{k=1}^{\infty}L^{k}(A)$$
where $A\in \mathbb{C}[[x]]$, $f=\sum_{k=0}^{n}{f_{k}\partial^k}+\sum_{k=1}^{\infty}{f_k\partial^{-k}}$ and 
$L(\sum_{i=0}^{\infty}c_ix^i)=\sum_{i=0}^{\infty}c_i\frac{x^{i+1}}{i+1}$, is well-defind and yields to the embedding.
\end{proof}
\end{subsection}
                                                      \begin{subsection}{Calkin co-cycle}
 It is easy to see that $I(F)=C\bigoplus [I(F),I(F)]$ where $[I(F),I(F)]$ is the linear space spanned by elements of the form $[L,L_1]=LL_1-L_1L$, $L,L_1\in I(F)$. 
 So the projection map on $C$ gives us a trace map on $I(F)$, denoted by $tr:I(F)\to C$. 
Consider $F(\partial)$ as a Lie algebra over $C$.
Assume that $\alpha:F(\partial)\to End_C(F)$ is a linear map such that $\pi  \alpha=\theta$, where $\pi :End_C(F)\to M_C(F)$ is the quotient map. 
Define $\sigma_{\alpha}$ on $F(\partial)$  as follows,\\
$$\sigma_{\alpha}(P,Q)=tr(\alpha([P,Q])-[\alpha(P),\alpha(Q)]).$$
It is easy to check that,
\begin{proposition}
(a) $\sigma_{\alpha}$ is a Lie algebra 2-cocycle on $F(\partial)$.\\
(b) The class of this 2-cocycle $[\sigma_{\alpha}]$ in $H^{2}_{Lie}(F(\partial))$ does not depend on $\alpha$.
\end{proposition}
\begin{definition}
The 2-cocycle constructed above is called the Calkin co-cycle. 
\end{definition}
 
 \end{subsection}
 
                                                       \begin{subsection}{Calkin and Kac co-cycles}

In this section we show the relation between the Calkin co-cycle and the Kac co-cycle.

We recall that the Lie algebra of differential operators on the circle, i.e. 
$\mathbb{C}[z,z^{-1}][\frac{\partial}{\partial z}]$, has a nontrivial central extension given by the Kac 2-cocycle( See [Kac]). 
We construct another 2 co-cycle coming from the Calkin co-cycle.\\
Consider differential ring $\mathbb{C}[x]$ whose derivative is $\frac{\partial}{\partial x}$. Since the derivative is onto, one can see that we have an embedding
$\theta:\mathbb{C}[x]\langle \partial,\partial^{-1}\rangle \to M_{\mathbb{C}}(\mathbb{C}[x])$ as the one in section 5.2( detailes are left to the reader). 
Therefore we obtain a 2
co-cycle on $\mathbb{C}[x]\langle \partial,\partial^{-1}\rangle$. In order to calculate this co-cycle, we define 
$\alpha_0:\mathbb{C}[x]\langle \partial,\partial^{-1}\rangle \to End_{\mathbb{C}}(\mathbb{C}[x])$ as follows. Every element of 
$\mathbb{C}[x]\langle \partial,\partial^{-1}\rangle $ can be written as $\sum_{i=0}^{n}{a_i\partial^i}+\sum_{j=1}^{m}b_j\partial^{-1}c_j$ where 
$a_i,b_j,c_j\in \mathbb{C}[x]$. We define 
$$\alpha_0(\sum_{i=0}^{n}{a_i\partial^i}+\sum_{j=1}^{m}b_j\partial^{-1}c_j)=
\sum_{i=0}^{n}{a_iD^i}+\sum_{j=1}^{m}{b_jLc_j}$$
where $L(x^k)=\frac{x^{k+1}}{k+1}$($k\geq 0$) and $D=\frac{\partial}{\partial x}$. It is easy to see that $\pi \alpha_0=\theta$. Set $\sigma_0=\sigma_{\alpha_0}$.
\begin{lemma} 
For any $a,a_1,b,b_1,c,c_1 \in  \mathbb{C}[x] $ and $ n,n_1\geq 0$,\\
(a)  $\sigma_0(a\partial^{n},a_1\partial^{n_1})=\sigma_0(b\partial^{-1}c,b_1\partial^{-1}c_1)=0$,\\
(b) $\sigma_0(a\partial^{n},b\partial^{-1}c)=\left ( \sum_{i=0}^{n-1}{(-1)^{i+1} (b(ac)^{(i)})^{(n-i-1)}} \right )(0)$
\end{lemma}
\begin{proof}
 (a) This is easy to see because $\alpha_0|_{\mathbb{C}[x]\langle \partial \rangle}$ and $\alpha_0|_{\mathbb{C}[x]\langle \partial^{-1} \rangle}$ 
 are in fact ring homomorphisms.\\
 (b) Setting $P=a\partial^{n}$ and $Q=b\partial^{-1}c$, one can easily verify that,
$$\alpha_0([P,Q])-[\alpha_0(P),\alpha_0(Q)]= \alpha_0(Q)\alpha_0(P)-\alpha_0(QP).$$
Suppse that $h:\mathbb{C}[x]\to \mathbb{C}$ is the evaluation map at zero. Then using  proposition 4.2, we have
 $$\alpha_0(Q)\alpha_0(P)-\alpha_0(QP)=
 \sum_{i=0}^{n-1}{(-1)^{i}b(ac)^{(i)}LD^{n-i}}-\sum_{i=0}^{n-1}{(-1)^{i}b(ac)^{(i)}D^{n-i-1}}=$$
$$\sum_{i=0}^{n-1}{(-1)^{i+1}b(ac)^{(i)}(1-LD)D^{n-i-1}}.$$
Since $h=1-LD$, we have,
$$\sigma_0(a\partial^{n},b\partial^{-1}c)=tr(\sum_{i=0}^{n-1}{(-1)^{i+1}b(ac)^{(i)}hD^{n-i-1}})$$
 $$=\sum_{i=0}^{n-1}{tr((-1)^{i+1}b(ac)^{(i)}hD^{n-i-1})}=\sum_{i=0}^{n-1}{(-1)^{i+1} h((b(ac)^{(i)})^{(n-i-1)})}$$
 $$=\left ( \sum_{i=0}^{n-1}{(-1)^{i+1} (b(ac)^{(i)})^{(n-i-1)}} \right )(0) .$$

 \end{proof}
One can see that there is a $\mathbb{C}$-algebra isomorphism 
 $$g:\mathbb{C}[z,z^{-1}][\frac{\partial}{\partial z}]\to \mathbb{C}[x]\langle \partial,\partial^{-1}\rangle$$ 
sending $z$ to $\partial$ and $\frac{\partial}{\partial z}$ to $-x$. This isomorphism gives a  2-cocycle on $\mathbb{C}[z,z^{-1}][\frac{\partial}{\partial z}]$ 
from the Calkin cocycle, denoted by $\sigma_1$. We need a lemma,
\begin{lemma}
For any $n<0$ we have,
$$\partial^{n}=\sum_{i=0}^{-n-1}{\frac{(-1)^{i}}{i!(-n-i-1)!}x^{-n-i-1}\partial^{-1}x^{i}}$$
\end{lemma}   
Using lemma 5.6. we have the following lemma,
\begin{lemma}For any $m,r\in \mathbb{Z}$ and $n,s\geq 0$,
 $\sigma_1(z^{m}\frac{\partial^{n}}{\partial z^{n}},z^{r}\frac{\partial^{s}}{\partial z^{s}})$ is zero if $rm\geq 0$ or $n+s\neq m+r$.
\end{lemma}
\begin{proof}
We have,
 $$\sigma_1(z^{m}\frac{\partial^{n}}{\partial z^{n}},z^{r}\frac{\partial^{s}}{\partial z^{s}})=(-1)^{n+s}
\sigma_0(\partial^{m}x^n,\partial^{r}x^s)$$
where $m,r\in \mathbb{Z}$ and $n,s\geq 0$. So it is zero unless $rm<0$. We may assume that $r<0$ and $m>0$. Using the above lemma, we have,
$$\sigma_1(z^{m}\frac{\partial^{n}}{\partial z^{n}},z^{r}\frac{\partial^{s}}{\partial z^{s}})=$$
$$(-1)^{n+s}\sigma_0(\sum_{i=0}^{m}{{m\choose i}(x^{n})^{(i)}\partial^{m-i}},\sum_{j=0}^{-r-1}{\frac{(-1)^{i}}{i!(-r-i-1)!}x^{-r-i-1}\partial^{-1}x^{i+s}})$$
$$=\sum_{i=0}^{m}\sum_{j=0}^{-r-1}{{m\choose i}\frac{(-1)^{i}}{i!(-r-i-1)!}\sigma((x^{n})^{(i)}\partial^{m-i},x^{-r-i-1}\partial^{-1}x^{i+s})}.$$
 Using lemma 5.4. and a simple calculation, we can see that $\sigma_1(z^{m}\frac{\partial^{n}}{\partial z^{n}},z^{r}\frac{\partial^{s}}{\partial z^{s}})$ 
is zero if $n+s\neq m+r$.
 \end{proof}
 Now we can explain the relation between the Calkin cocycle and Kac cocycle. First some lemmas,
\begin{lemma}
If $L$ is  a Lie algebra over $\mathbb{C}$ and $\beta$ a nonzero 2-cocycle on $L$ such that $\{[a,b]|\beta(a,b)=0\}$ generates $L$ as  
a $\mathbb{C}$-vector space, then $\beta$ is nonzero in $H^{2}(L)$.
 \end{lemma}
 See [Li] Lemma 1.
 \begin{lemma}
 We have $[\mathbb{C}[x]\langle \partial \rangle,\mathbb{C}[x]\langle \partial \rangle]=\mathbb{C}[x]\langle \partial \rangle$, 
 $[\mathbb{C}[x]\langle \partial^{-1}\rangle,\mathbb{C}[x]\langle  \partial^{-1}\rangle]=\mathbb{C}[x]\langle  \partial^{-1}\rangle\partial^{-2}$ and 
 $[a\partial,\partial^{-1}]=\partial^{-1}a'$ 
for any $a\in \mathbb{C}[x]$.
 \end{lemma}
 \begin{proof}
 A simple calculation.
 \end{proof}
 \begin{lemma}
 $H^{2}_{Lie}(\mathbb{C}[z,z^{-1}][\frac{\partial}{\partial z}])$ is one dimensional over $\mathbb{C}$.
 \end{lemma}
  See [Gelfand] for a proof.
  
 Combining these lemmas we have,
 \begin{proposition}
 The cocycle $\sigma_1$ is nontrivial.
 \end{proposition}
 \begin{lemma}
 Lemma 5.4. shows that $\sigma_0$ is zero on $\mathbb{C}[x]\langle \partial \rangle$ and $\mathbb{C}[x]\langle  \partial^{-1}\rangle$. Furthermore
$$\sigma_0(a\partial^{n},\partial^{-1})=a(0)$$
for any $a\in \mathbb{C}[x]$. Now lemma 5.9 together with lemma 5.10. gives the conclusion.
 \end{lemma}
  There is another way to define the unique one-dimensional central extension of $\mathbb{C}[z,z^{-1}][\frac{\partial}{\partial z}]$ given by Kac cocycle $\tau$
  ( see [Kac] theorem 5.3.). Then,
  \begin{proposition}
  We have $\sigma_1=\tau$.
  \end{proposition}
  \begin{proof}
  By lemma 5.11., $[\sigma_1]=c[\tau]$ in $H^{2}_{Lie}(\mathbb{C}[z,z^{-1}][\frac{\partial}{\partial z}])$ for some $c\in \mathbb{C}$. But this implies that $\tau=c\sigma_1$, 
  becuase of lemma 5.9. Finally a simple calculation shows that $c=1$.
  \end{proof}
   
\end{subsection}

\end{section}

\textbf{e-mail: masood.aryapoor@yale.edu}

 \end{document}